\documentclass{IEEEtran4PSCC}

\ifCLASSINFOpdf
   \usepackage[pdftex]{graphicx}
\else
   \usepackage[dvips]{graphicx}
\fi
\usepackage[cmex10]{amsmath}
\usepackage{amsfonts,multirow,array,balance}
\newcolumntype{P}[1]{>{\centering\arraybackslash}p{#1}}

\usepackage{xparse}
\newsavebox{\fminipagebox}
\NewDocumentEnvironment{fminipage}{m O{\fboxsep}}
 {\par\kern#2\noindent 
  \vspace{-3ex}
 \begin{center}\begin{lrbox}{\fminipagebox}
  \begin{minipage}{#1}\ignorespaces}
 {\end{minipage}\end{lrbox}%
  \makebox[#1]{%
    \kern\dimexpr-\fboxsep-\fboxrule\relax
    \fbox{\usebox{\fminipagebox}}%
    \kern\dimexpr-\fboxsep-\fboxrule\relax
  }\par\kern#2 \end{center}
  \vspace{-1ex}
 }

\newcounter{ProblemCounter}
\renewcommand\theProblemCounter{\arabic{ProblemCounter}}

\DeclareMathOperator{\diag}{dg}

\DeclareMathOperator{\trace}{Tr}

\DeclareMathOperator{\real}{Re}
\DeclareMathOperator{\imag}{Im}

\newtheorem{assumption}{Assumption}

\newtheorem{lemma}{Lemma}

\newcommand \bzero{\mathbf{0}}
\newcommand \bone{\mathbf{1}}

\newcommand \be{\mathbf{e}}

\newcommand \bi{\mathbf{i}}

\newcommand \bp{\mathbf{p}}

\newcommand \bu{\mathbf{u}}
\newcommand \bv{\mathbf{v}}

\newcommand \bx{\mathbf{x}}
\newcommand \by{\mathbf{y}}

\newcommand \bD{\mathbf{D}}
\newcommand \bE{\mathbf{E}}

\newcommand \bH{\mathbf{H}}
\newcommand \bI{\mathbf{I}}

\newcommand \bL{\mathbf{L}}
\newcommand \bM{\mathbf{M}}

\newcommand \bU{\mathbf{U}}
\newcommand \bV{\mathbf{V}}
\newcommand \bW{\mathbf{W}}
\newcommand \bX{\mathbf{X}}
\newcommand \bY{\mathbf{Y}}
\newcommand \bZ{\mathbf{Z}}

\newcommand \bdelta{\boldsymbol{\delta}}

\newcommand \bGamma{\mathbf{\Gamma}}

\newcommand \bLambda{\mathbf{\Lambda}}

\newcommand \mcE{\mathcal{E}}

\newcommand \mcG{\mathcal{G}}

\newcommand \mcL{\mathcal{L}}

\newcommand \mcN{\mathcal{N}}

\newcommand \mcS{\mathcal{S}}

\newcommand \bbp{\bar{\mathbf{p}}}

\newcommand \bbdelta{\bar{\boldsymbol{\delta}}}

\renewcommand{\d}[1]{\ensuremath{\operatorname{d}\!{#1}}}

\makeatletter
\let\old@ps@headings\ps@headings
\let\old@ps@IEEEtitlepagestyle\ps@IEEEtitlepagestyle
\def\psccfooter#1{%
    \def\ps@headings{%
        \old@ps@headings%
        \def\@oddfoot{\strut\hfill#1\hfill\strut}%
        \def\@evenfoot{\strut\hfill#1\hfill\strut}%
    }%
    \def\ps@IEEEtitlepagestyle{%
        \old@ps@IEEEtitlepagestyle%
        \def\@oddfoot{\strut\hfill#1\hfill\strut}%
        \def\@evenfoot{\strut\hfill#1\hfill\strut}%
    }%
    \ps@headings%
}
\makeatother

\begin{document}

\title{Optimal Power Flow Schedules with Reduced Low-Frequency Oscillations}

\author{
\IEEEauthorblockN{Manish K. Singh and Vassilis Kekatos}
\IEEEauthorblockA{Bradley Dept. of Electrical \& Computer Engnr.\\ Virginia Tech, Blacksburg, VA 24061, USA\\
\{manishks, kekatos\}@vt.edu}
}

\maketitle
\begin{abstract}
The dynamic response of power grids to small events or persistent stochastic disturbances influences their stable operation. Low-frequency inter-area oscillations are of particular concern due to insufficient damping. This paper studies the effect of the operating point on the linear time-invariant dynamics of power networks. A pertinent metric based on the frequency response of grid dynamics is proposed to quantify power system's stability against inter-area oscillations. We further put forth an optimal power flow formulation to yield a grid dispatch that optimizes this novel stability metric. A semidefinite program (SDP) relaxation is employed to yield a computationally tractable convex problem. Numerical tests on the IEEE-39 bus system demonstrate that the SDP relaxation is exact yielding a rank-1 solution. The relative trade-off of the proposed small-signal stability metric versus the generation cost is also studied. 
\end{abstract}


\section{Introduction}
\allowdisplaybreaks
Power systems experience fluctuations of small magnitudes at all times, and are often exposed to large disturbances as well. Stability, being the ability to withstand such disturbances, is naturally an extremely sought attribute. Power system stability studies have been classified in diverse ways depending upon the variables of interest, the magnitude of disturbances, the considered timeframe, as well as the order and accuracy of component modeling~\cite{stability_definition_1}, \cite{stability_definition_2}, \cite{kundur}, \cite{Chow2020Book}. In large power transmission networks, rotor-angle stability has been widely investigated, where transient and small-signal stability studies explore the effect of large and small disturbances, respectively.

Transient stability studies are primarily conducted via time-series simulations or direct methods based on energy-type Lyapunov functions that characterize stability regions~\cite{Chow2020Book}, \cite{Turitsyn2016lyapunov}. Such approaches may not be practically scalable, and thus oftentimes confine transient stability analyses to a limited set of probable contingencies~\cite{Zimmerman2000stability}. Transient stability is known to be strongly related to the operating point, and is actually studied accordingly. Attempts have also been made at formulating OPF problems that yield operating points with guaranteed transient stability~\cite{Zimmerman2000stability}, \cite{chen2020robust}.

For small-signal stability, the power network is oftentimes approximated by a linear time-invariant (LTI) system obtained upon linearizing the chosen dynamical model and the power flow equations at a given operating point. Thereafter, the system's small-signal stability is assessed by investigating the eigenvalues of the related \emph{state matrix}; specifically ensuring that they possess negative real parts. Thus, the stability at an operating point is often ensured by minimizing the \emph{spectral abscissa}, that is the real-part of the eigenvalue closest to the imaginary axis~\cite{kundur2004generation}, \cite{conejo2011smallsignal}. Instead of obtaining operating points explicitly, several works incorporate Lyapunov function-based approaches to design stability-ensuring control mechanisms~\cite{Taha2019Linfinity}, \cite{Baros2021distributed}. Alternatively, reference~\cite{Chakrabortty21} seeks a stability-promoting OPF dispatch by adjusting the reactive power setpoints of loads. Although the approach considers detailed power system dynamics, it engages only a subset of loads and are adjusted by gradient updates that may miss global optimality.

Within the purview of small-signal stability, there has been a special focus on limiting inter-area oscillations, which often involve groups of synchronous machines oscillating against each other. Inter-area oscillations are known to occupy the lower frequency spectrum of grid signals~\cite{kundur}. Incidents of sustained low-frequency oscillations, oftentimes of unknown origin, are frequently reported in large power systems~\cite{Chow2000damping}. If not properly controlled, inter-area oscillations can trigger catastrophic outages such as the 1996 West Coast Blackout in the United States~\cite{Kosterev99blackout}. While intra-area oscillations can be damped effectively by local controllers, inter-area oscillations are harder to control, thus exacerbating the threat~\cite{PSS92}.

Early approaches for limiting inter-area oscillations were centered at designing power system stabilizers with local control signals. While the early works had limited efficacy, tremendous advancements have led to the rich field of \emph{wide-area damping control}. The state of the art involves meticulous signal selection, control design and placement, and diverse actuators such as FACTS devices, energy storage, wind generators, load control, and control of high-voltage direct current lines; see~\cite{Jianming20load} and references therein. A typical control design approach undertaken in the aforementioned works involves linearizing the dynamical system at a \emph{fixed} operating point; selecting the concerning mode(s); and designing control parameters to alter the related closed-loop poles as desired. Despite ignoring the updates in operating point, the resulting controls have been traditionally effective because the dominant inter-area modes were known to not change significantly~\cite{Chow2000damping}. With increasing agility in power grid operation however, recent works advocate re-evaluating the linearized models based on frequent grid updates or online system identification~\cite{Jianming20load}, \cite{Sudipta20pss}.

In this work, we further emphasize on the dependence of power system dynamic behavior on the operating point and seek to exploit this dependence in improving small-signal stability. Indeed, the works aiming at determining operating points that minimize spectral abscissa capitalize on the aforementioned dependence, albeit at high computational cost~\cite{conejo2011smallsignal}. This work focuses on determining operating points that limit inter-area oscillations. Building on the observation that inter-area oscillations occupy the lower end of system frequency response, we propose a novel optimal power flow (OPF) formulation that minimizes the energy of low-frequency modes. The proposed indirect approach of identifying inter-area modes via frequency response instead of the mode-shapes obviates the need for cumbersome computation of eigenvectors. Furthermore, the proposed approach is amenable to a semidefinite program (SDP) relaxation, and integrates conveniently with traditional OPFs that minimize generation cost. Thus, we used the resulting formulation for trade-off analysis between small-signal stability and generation cost. Numerical tests on the IEEE 39-bus system interestingly yielded nearly rank-1 minimizers for all instances, rendering the SDP relaxation exact.

The rest of this work is organized as follows. Section~\ref{sec:model} reviews linearized swing dynamics. Section~\ref{sec:oscillations} studies the frequency-domain characteristics of swing dynamics, quantifies the energy of inter-area oscillations under ambient dynamics, and provides an SDP representation for the new metric. Section~\ref{sec:opf} puts forth an SDP-based OPF to reduce the energy of inter-area oscillations by adjusting the operating point. The exactness of the relaxation and the trade-off between generation cost and stability are evaluated through extensive numerical tests in Section~\ref{sec:tests}.

\emph{Notation}: lower- (upper-) case boldface letters denote column vectors (matrices). Calligraphic symbols are reserved for sets. Symbol $^{\top}$ ($^{H}$) stands for (complex) transposition, vectors $\mathbf{0}$ and $\mathbf{1}$ are the all-zeros and all-ones vectors or matrices.

\section{Power System Modeling}\label{sec:model}
This section reviews swing dynamics, relates internal to external voltages, and linearizes swing equations to establish the link between dynamics and the steady-state operating point of the power system. A power transmission network with $N$ buses can be represented by an undirected connected graph $\mcG=(\mcN,\mcL)$, whose nodes $n\in{\mcN}:=\{1,\ldots,N\}$ correspond to buses, and edges $\ell=(m,n)\in\mcL$ to transmission lines. 	Set $\mcN$ is a collection of generation, load, and zero-injection buses. In high-voltage transmission systems, generator buses primarily host synchronous machines. Load buses however oftentimes denote the connection to sub-transmission networks, which indeed host numerous medium-sized synchronous generators and motors, alongside other types of loads. We therefore collectively refer to generation and load buses as synchronous buses, collect them in $\mcS\subset\mcN$, and index them by $n=1,\dots,S$. Non-synchronous buses form set $\bar{\mcS}=\mcN\setminus\mcS$ and are indexed by $n=S+1,\dots,N$. 	
	
We first review the angular dynamics associated with the buses in $\mcS$. Let the synchronous machine connected to bus $n\in\mcS$ be associated with a rotor speed $\dot{\delta}_n$, inertia constant $M_n>0$, and damping coefficient~$D_n>0$; see~\cite{kundur} for details. The rotor speed $\dot{\delta}_n$ is the deviation from the nominal, which ideally should be close to zero. Swing dynamics for the synchronous machine connected to bus~$n$ predicate~\cite{kundur}
\begin{equation}\label{eq:swing}
M_n\ddot{\delta}_n + D_n\dot{\delta}_n= p_n^\text{in}-p_n 
\end{equation} 
where $p_n^\text{in}:=p_n^g-p_n^d$ is the difference between the mechanical power input $p_n^g$ and the local power demand $p_n^d$ at bus $n$, while $p_n$ is the electric power flowing from bus $n$ into the grid. It is worth stressing that swing dynamics relate to the \emph{internal} voltages of synchronous machines denoted by $e_n = E_n\angle \delta_n$ for all $n\in\mcS$. When studying dynamics, power injections $p_n$ are typically expressed in terms of internal voltages as well.

On the contrary, OPF formulations naturally represent synchronous buses using their \emph{external voltages} denoted here by $v_n = V_n\angle \theta_n$. Voltage phasors $v_n$ are defined for all $n\in\mcN$, that is for synchronous and zero-injection buses alike. In OPF studies, power injections $p_n$ and constraints on line flows or voltage magnitudes are all expressed as functions of $v_n$'s.

Let $\bv$ and $\bi$ be the $N$-length vectors collecting respectively the complex external voltages and injection currents across all buses $n\in\mcN$. Vector $\be$ collects the complex internal voltages at all buses $n\in\mcS$. Given line impedances and shunt susceptances, one can derive the $N\times N$ bus admittance matrix $\bY$, satisfying Ohm's law $\bi=\bY\bv$. 

Let us block partition Ohm's law on buses belonging to $\mcS$ and $\bar{\mcS}$. Before doing so, note that each internal bus is connected to its associated external bus $n$ via reactance $x_n>0$. If all such reactances are collected  by a diagonal matrix as $\bY_{\mcS}:=\diag(\{1/(jx_n)\})$, then Ohm's law can be partitioned as~\cite{Aranya2018graph}
	\begin{equation}\label{eq:ohm}
	    \begin{bmatrix}
	\bY_{\mcS,\mcS} & \bY_{\mcS,\bar{\mcS}} \\
	\bY_{\mcS,\bar{\mcS}}^\top &\bY_{\bar{\mcS},\bar{\mcS}}
	\end{bmatrix}\begin{bmatrix}
	    \bv_{\mcS}\\
	    \bv_{\bar{\mcS}}
	\end{bmatrix}=\begin{bmatrix}
	    \bY_{\mcS}(\be-\bv_{\mcS})\\
	    \bzero
	\end{bmatrix}.
	\end{equation}
Eliminating $\bv_{\bar{\mcS}}$ via the so termed \emph{Kron reduction} yields
	\begin{equation}\label{eq:ve}
	    \bv_{\mcS}=\bGamma\bY_{\mcS}\be
	\end{equation}
where matrix $\bGamma:=(\bY_{\mcS,\mcS}+\bY_{\mcS}-\bY_{\mcS,\bar{\mcS}}\bY_{\bar{\mcS},\bar{\mcS}}^{-1}\bY_{\mcS,\bar{\mcS}}^\top)^{-1}$ is known to exist~\cite{Dorfler13}. According to~\eqref{eq:ve}, the external voltages at synchronous buses are linear functions of internal voltages.
	
The active power injection $p_n$ from bus $n$ into the network can be interpreted as the active flow over reactance $x_n$ between the internal voltage $e_n$ and the external voltage $v_n$ as
    \begin{equation}\label{eq:p-ev}
        p_n=\frac{E_nV_n}{x_n}\sin(\delta_n-\theta_n)=\frac{1}{x_n}\imag\{e_nv_n^*\}.
    \end{equation}
    Eliminating $v_n$ from \eqref{eq:ve} in \eqref{eq:p-ev}, one obtains
    \begin{align}\label{eq:p-ee1}
        p_n=\sum_{m=1}^S\frac{E_nE_m}{x_nx_m}\big[&\real\{\Gamma_{nm}\}\cos(\delta_n-\delta_m)\nonumber\\&-\imag\{\Gamma_{nm}\}\sin(\delta_n-\delta_m)\big].
    \end{align}
Owing to the observation that $|\real(\Gamma_{nm})|\ll|\imag(\Gamma_{nm})|$, the assumption of a lossless transmission system is oftentimes invoked to approximate
        \begin{align}\label{eq:p-ee2}
        p_n=\sum_{m=1}^S\frac{E_nE_m}{\gamma_{nm}}\sin(\delta_n-\delta_m)
    \end{align}
where $\gamma_{nm}:=-(x_nx_m)/\imag(\Gamma_{nm})>0$ serves as the \emph{effective reactance} between the synchronous machines connected to buses $m$ and $n$~\cite{Aranya2018graph}. By symmetry, it holds that $\gamma_{nm}=\gamma_{mn}$. Although~\eqref{eq:swing} provides individual machine dynamics, the power injection $p_n$ couples the system of $S$ machines due to~\eqref{eq:p-ee2}. 

Let us linearize the swing dynamics of~\eqref{eq:swing} at the equilibrium. To this end, assume constant voltage magnitudes and consider small perturbations from the steady-state operating schedule $\bar{p}_n^\text{in}$ to $\bar{p}_n^\text{in}+\Delta p_n^\text{in}$. Such changes induce in turn perturbations from $(\bar{p}_n,\bar{\delta}_n)$ with $\bar{p}_n=\bar{p}_n^\text{in}$ to $(\bar{p}_n+\Delta p_n,\bar{\delta}_n+\Delta\delta_n)$. Linearizing the right-hand side of~\eqref{eq:p-ee2} and collecting variables across all synchronous buses in vectors gives
	\begin{equation}\label{eq:linPF}
	\bbp^\text{in}+\Delta\bp^\text{in}-\bbp-\Delta\bp=\Delta\bp^\text{in}-\bL_{\bbdelta}\Delta\bdelta 
	\end{equation}
since $\bbp^\text{in}=\bbp$ at equilibrium, and where $\bL_{\bbdelta}$ is the Jacobian matrix of \eqref{eq:p-ee2} with respect to the internal voltage angles $\bdelta$ evaluated at $\bbdelta$. Its entries are
\begin{equation}\label{eq:Laplacian}
[\bL_{\bbdelta}]_{n,m}=\left\{\begin{array}{ll}
\sum_{k\neq n} \frac{E_nE_k}{\gamma_{nk}}\cos(\bar{\delta}_n-\bar{\delta}_k)&,~m=n\\
\\
-\frac{E_nE_m}{\gamma_{nm}}\cos(\bar{\delta}_n-\bar{\delta}_m)&,~ m\neq n.
\end{array}\right.
\end{equation}
Heed that $\gamma_{nm}>0$, and for admissible power system operating points, it holds that $|\bar{\delta}_n-\bar{\delta}_m|<\pi/2$, implying $\cos(\bar{\delta}_n-\bar{\delta}_m)>0$. Therefore, the Jacobian $\bL_{\bbdelta}$ is in fact a Laplacian matrix with positive edge weights given by $E_nE_m\gamma_{nm}^{-1}\cos(\bar{\delta}_n-\bar{\delta}_m)$ for edge $(n,m)\in\mcL$. With the aforementioned observation, matrix $\bL_{\bbdelta}$ can be shown to be positive semidefinite~\cite{Low18}.

Substituting \eqref{eq:linPF} into \eqref{eq:swing} yields
	\begin{equation}\label{eq:swing2}
	\bM\ddot{\bdelta} + \bD\dot{\bdelta}+\bL_{\bbdelta}\Delta\bdelta=\Delta\bp^\text{in}
	\end{equation}
where diagonal matrices $\bM$ and $\bD$ carrying the inertia and damping coefficients on their main diagonals; and $\ddot{\bdelta}$ and $\dot{\bdelta}$ collect rotor accelerations (rate of change of frequencies, ROCOF) and speeds across buses. Note the time derivatives of $\bdelta$ and $\Delta\bdelta$ coincide. Therefore, the LTI system in~\eqref{eq:swing2} involves the deviations of states and inputs from their nominal values. 

For notational brevity, we henceforth omit $\Delta$'s; simplify $\bp^\text{in}$ as $\bp$; and drop the dependence of $\bL_{\bbdelta}$ on $\bbdelta$ to get
\begin{equation}\label{eq:swing3}
\bM\ddot{\bdelta} + \bD\dot{\bdelta}+\bL\bdelta= \bp.
\end{equation}
As evident by \eqref{eq:swing3}, the steady-state operating point affects swing dynamics through the Laplacian matrix $\bL$. Matrix $\bL$ is evaluated at the equilibrium values of internal voltages. Heed that it depends not only on the angles $\delta_n$'s, but also the magnitudes $E_n$'s of $\be$. This matrix will later serve as the link between swing oscillations and the OPF schedule.

\section{Metrics for Inter-area Oscillations}\label{sec:oscillations}
This section decouples the swing dynamics of \eqref{eq:swing3} and studies their frequency-domain characteristics to quantify the energy of inter-area oscillations under ambient dynamics. This metric will be incorporated later into an OPF to yield stability-informed schedules.

\subsection{Decoupled Swing Dynamics}
Let us study the dynamics of~\eqref{eq:swing3} in the frequency domain. This not only facilitates analytical convenience, but also allows us to explicitly focus on low-frequency inter-area oscillations. The transfer function of the LTI system in \eqref{eq:swing3} is
\begin{equation}\label{eq:tf1}
\bH(s)=\left(s^2\bM+s\bD+\bL\right)^{-1}
\end{equation}
with $s$ being the complex frequency of the Laplace domain. This transfer function simplifies significantly under the next assumption, which is adopted frequently to approximate power system dynamics~\cite{Low18},~\cite{Paganini19}.
	
\begin{assumption}\label{as:1}
The ratio of each generator's damping coefficient to its inertia is constant or $\bD=\gamma \bM$ for some $\gamma >0$.
\end{assumption}

This assumption relies on the fact that both inertia and damping coefficients of a synchronous machine scale with the machine's power rating~\cite{Poolla2017inertia}. Under this assumption, the transfer function of swing dynamics can be rewritten as~\cite{Paganini19}
\[\bH(s)=\bM^{-1/2}\bU\left(s^2\bI+s\gamma\bI+\bLambda\right)^{-1}\bU^\top\bM^{-1/2}\]
where $\bU \bLambda \bU^{\top}$ is the eigenvalue decomposition of the positive semidefinite matrix \begin{equation}\label{eq:L-LM}
    \bL_M := \bM^{-1/2} \bL \bM^{-1/2}.
\end{equation} The diagonal matrix $\bLambda$ carries the eigenvalues of $\bL_M$ sorted in increasing order as $0=\lambda_1<\lambda_2\leq \ldots\leq \lambda_S$, and the orthonormal matrix $\bU$ carries the associated eigenvectors as its columns. 

Let us now transform the original inputs/states of \eqref{eq:swing} to the \emph{eigeninputs/eigenstates}~\cite{Zhu18},~\cite{Paganini19}
\begin{equation}\label{eq:trans}
\bx:=\bU^{\top} \bM^{-1/2} \bp~~~ \text{and}~~~ \by := \bU^{\top} \bM^{1/2} \bdelta.
\end{equation}
Then, the swing dynamics of \eqref{eq:swing} transform to
\begin{align}\label{eq:eigenspace}
\ddot{\by} + \gamma \dot{\by} + \bLambda \by = \bx.
\end{align}
Because $\bLambda$ is diagonal, the MIMO swing dynamics of \eqref{eq:swing3} decouple into $S$ independent SISO \emph{eigensystems}. Eigensystem $i$ is described by the differential equation
\begin{equation}\label{eq:diag}
\ddot{y}_i + \gamma \dot{y}_i + \lambda_i y_i = x_i.
\end{equation}
Its impulse response can be computed as~\cite{Zhu21}
\begin{equation}\label{eq:h}
    h_i(t)=\frac{1}{r_i}\left( e^{c_i t} -e^{d_i t}\right) u(t)
\end{equation}
where $u(t)$ is the unit step function and \[r_i=\sqrt{\gamma^2-4\lambda_i},\quad c_i=\frac{-\gamma+r_i}{2},\quad d_i=\frac{-\gamma-r_i}{2}.\]
Note that $r_i$ is imaginary when $4\lambda_i>\gamma^2$, yet $h_i(t)$ in \eqref{eq:h} remains real-valued.

Focus on eigensystem $i=1$ corresponding to eigenvalue $\lambda_1=0$. This system is only marginally stable. For this reason, it will be excluded from future discussions. This is without loss of generality as the contribution of this first eigensystem is inconsequential to swing dynamics. To explain this, note that the eigenvector associated with eigenvalue $\lambda_1$ is
\begin{equation}\label{eq:u1}
\bu_1=\alpha_1\bM^{1/2}\bone~~~\text{where}~~~\alpha_1:=\left(\bone^\top\bM\bone\right)^{-1/2}.
\end{equation}
This is easy to verify as $\bL_M \bu_1=\bzero$ since $\bL\bone=\bzero$. Scaling by $\alpha_1$ ensures $\|\bu_1\|_2=1$. From \eqref{eq:trans}, we get that $\bdelta=\bM^{-1/2}\bU\by$, and so the contribution of eigenstate $y_1(t)$ to all nodal angles $\bdelta(t)$ is vector $\bM^{-1/2}\bu_1 y_1(t)=\alpha_1\bone y_1(t)$, which is a constant shift across all nodal angles. Because swing dynamics are shift-invariant as $\bL\bone=\bzero$, we can safely ignore $y_1(t)$.

We would like to focus on a reduced number of eigenstates associated with \emph{inter-area oscillations}~\cite{kundur}. Inter-area oscillations occupy the lower frequency spectrum of dynamic grid signals $\bdelta$, $\dot{\bdelta}$, and $\ddot{\bdelta}$, and can be observed over larger geographical areas or even across the entire power system~\cite[Ch.~12]{kundur}. Different from intra-area oscillations that can be damped effectively by local controllers, inter-area oscillations are harder to control and are thus of particular interest~\cite{PSS92}. 

\begin{figure}[thb]
    \centering
	\includegraphics[width=1\linewidth]{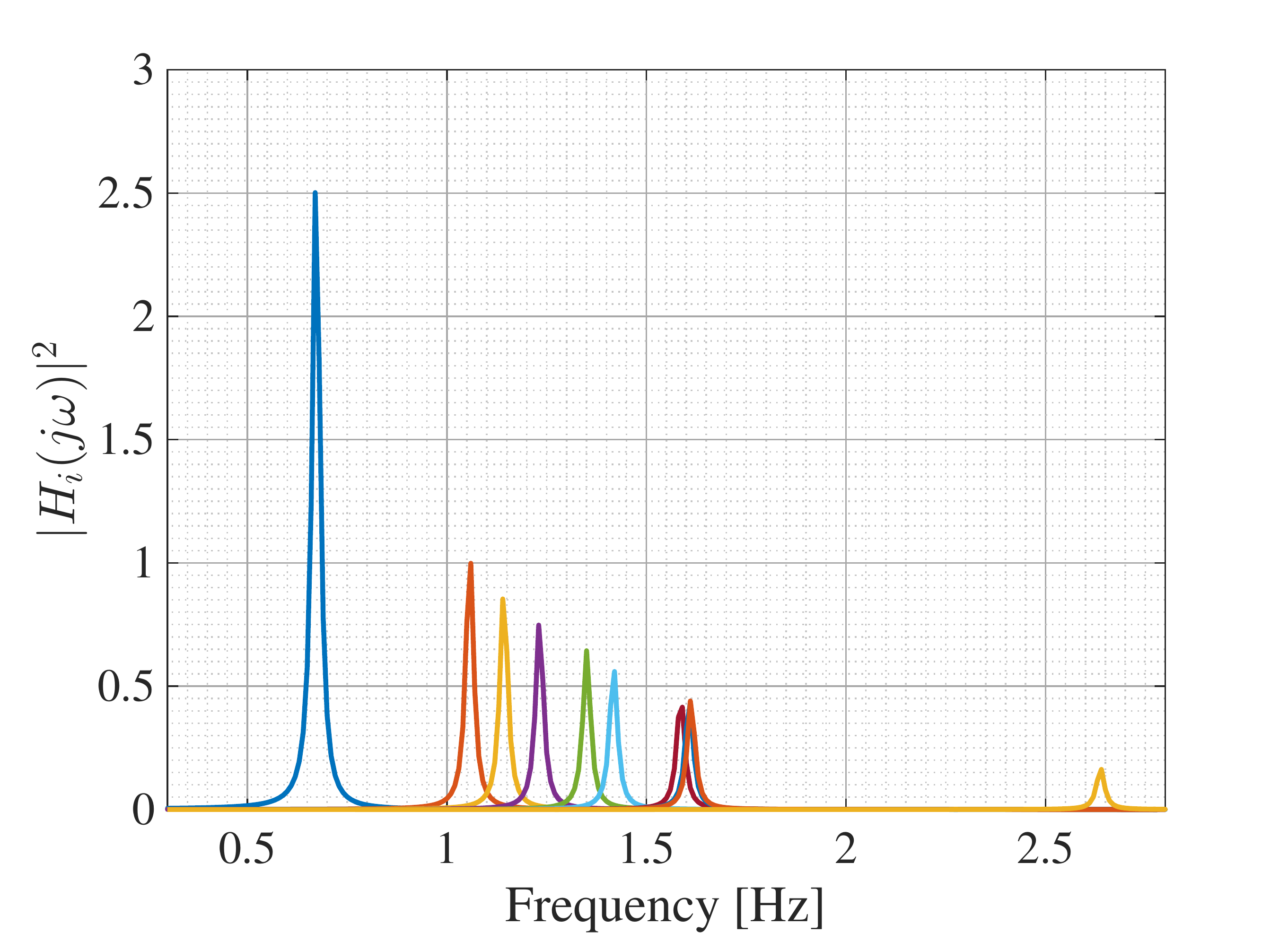}
	\caption{Frequency response magnitudes for the first ten eigensystems (excluding the first marginally stable one) of the IEEE 39-bus system.}
	\label{fig:H}       
\end{figure}

To target inter-area oscillations, we exploit key frequency-domain properties of eigenstates. In particular, the frequency response of the $i$-th eigensystem can be computed from \eqref{eq:h} for $i=2,\ldots,S$ as
\begin{equation}\label{eq:FR}
    |H_i(j\omega)|^2=\frac{1}{(\lambda_i-\omega^2)^2+\gamma^2 \omega^2}.
\end{equation}
Figure~\ref{fig:H} plots $|{H}_i(\omega)|^2$ for the eigensystems associated with the ten smallest eigenvalues of matrix $\bL_M$ for the IEEE 39-bus system, Kron reduced to a 29-bus system by eliminating 10 zero-injection buses. From \eqref{eq:FR}, one can observe that each eigensystem $i$ exhibits frequency selective behavior with its passband centered around the resonant frequency $\omega_i=\sqrt{\lambda_i-\frac{\gamma^2}{2}}$ at which the frequency response magnitude reaches its peak
\[|H_i(j\omega_i)|^2=\frac{1}{\gamma^2\left(\lambda_i-\frac{\gamma^2}{4}\right)}.\]
The aforementioned observation is graphically captured in Fig.~\ref{fig:H}, as well. Because practically $\lambda_i\gg \gamma^2/2$ for all $i\geq 2$,\footnote{In our numerical tests, we observed ${\lambda_2}/{\gamma^2}\simeq 829$.} resonant frequencies can be approximated as $\omega_i\simeq \sqrt{\lambda_i}$. At the resonant frequency $\omega_i$, the frequency response of the $i$-th eigensystem attains the value $|H_i(j\omega_i)|\simeq 1/(\gamma\sqrt{\lambda_i})$. The bandwidth of each eigensystem can be shown to scale approximately with $\gamma$. In other words, the passband for eigensystem $i$ shifts to higher frequencies with increasing $i$, whereas the frequency response gain decreases with increasing $i$. Due to this frequency selective behavior, the frequency content of eigenstate ${y}_i(t)$ can be approximately confined within the frequency band $\sqrt{\lambda_i}\pm \gamma$. Thus, to study inter-area oscillations, the operator can focus only on the $K$ eigenstates falling within the low-frequency band of interest. The eigensystems corresponding to resonant frequencies in the band of interest denoted by $[\underline{\omega},\overline{\omega}]$ comprise set $\mcE$, that is $\sqrt{\lambda_i}\in[\underline{\omega},\overline{\omega}]$ for all $i\in\mcE$. The remaining eigensystems constitute the complementary set $\bar{\mcE}=\mcN\setminus \mcE$. Heed that eigensystem $i=1$ with $\lambda_1=0$ belongs to $\bar{\mcE}$.

Given this grouping of eigensystems, let us collect the eigenvector/eigenvalue pairs corresponding to $\mcE$ in matrices $(\bU_\mcE,\bLambda_\mcE)$; and the remaining pairs in $(\bU_{\bar{\mcE}},\bLambda_{\bar{\mcE}})$, so that 
\[\bL_M = \bU_\mcE\bLambda_\mcE\bU_\mcE^\top + \bU_{\bar{\mcE}}\bLambda_{\bar{\mcE}}\bU_{\bar{\mcE}}^\top.\]
We can now express the original system states as
\begin{equation}\label{eq:dec}
\bdelta = \underbrace{\bM^{-1/2}\bU_\mcE\by_{\mcE}}_{:=\bdelta_\mcE} + \underbrace{\bM^{-1/2}\bU_{\bar{\mcE}}\by_{\bar{\mcE}}}_{:=\bdelta_{\bar{\mcE}}}
\end{equation}
where $\by_{\mcE}$ and $\by_{\bar{\mcE}}$ gather the entries of $\by$ related to the eigensystems in $\mcE$ and $\bar{\mcE}$, respectively. Recall that vectors $\by_{\mcE}$ and $\by_{\bar{\mcE}}$ carry signals with frequency content at different bands. Each angle $\delta_n(t)$ is a linear function of $\by_{\mcE}(t)$ and $\by_{\bar{\mcE}}(t)$. In essence, vector $\bdelta_\mcE$ can be obtained by filtering $\bdelta$ to keep frequencies falling only within the band of interest. To reduce the frequency content of $\bdelta$ across this band, the goal is to suppress signal $\by_{\mcE}$. Upon premultiplying \eqref{eq:dec} by $\bU_\mcE^\top\bM^{1/2}$ and $\bU_{\bar{\mcE}}^\top\bM^{1/2}$, and due to the orthonormality of $\bU$, it is easy to verify that
\[\by_\mcE=\bU_\mcE^\top\bM^{1/2}\bdelta~~~\text{and}~~~\by_{\bar{\mcE}}=\bU_{\bar{\mcE}}^\top\bM^{1/2}\bdelta~.\]
We next measure the energy of $\by_\mcE$ under a meaningful statistical model for ambient dynamics.

\subsection{A Metric of Inter-area Oscillations} 
Ambient dynamics are observed when the power network is driven by small random fluctuations in power injections $\bp(t)$ of \eqref{eq:swing3}. Lacking any particular statistical assumptions for such fluctuations, we adopt a non-informative prior and model $\bp(t)$ as a zero-mean white random process\footnote{The system response to a white random process features interesting relations to that obtained for node-wise impulse inputs; see~\cite{Poolla2017inertia} for details.} with covariance $\mathbb{E}[\bp(t+\tau)\bp^\top(t)]=\alpha\bM\delta_D(\tau)$ for some $\alpha>0$, where $\delta_D(\tau)$ is the Dirac delta function. The assumption here is that the variance of nodal power injection fluctuations scales with the per-bus power rating~\cite{Zhu18},~\cite{Zhu21}. Using \eqref{eq:trans}, the random vector of eigeninputs $\bx(t)$ is also zero-mean for all $t$. Its covariance matrix can be easily found to be~\cite{Paganini19},~\cite{GPCDC}
\begin{equation}\label{eq:xcov2}
\mathbb{E}[\bx(t+\tau)\bx^\top(t)]=\alpha\bI_N\delta_D(\tau).
\end{equation}
Without loss of generality, let us ignore $\alpha$ by setting it to unity. Equation~\eqref{eq:xcov2} implies that eigeninputs are white across time, and uncorrelated across eigenstates. If the SISO eigensystems in \eqref{eq:h} are driven by zero-mean random processes $x_i(t)$ uncorrelated with each other, their outputs $y_i(t)$ will be zero-mean random processes uncorrelated with each other as well. Mathematically, matrix $\mathbb{E}[\by(t+\tau)\by^\top(t)]$ is diagonal. Its diagonal entries can be computed as follows~\cite{Zhu21}, \cite{Jalali21}. If an LTI system [such as the $i$-th eigensystem in \eqref{eq:diag}] is driven by a wide-sense stationary (WSS) random process $x_i(t)$, its output is a WSS random process with covariance\footnote{A random process $z(t)$ is WSS if its mean $\mathbb{E}[z(t)]$ and covariance $\mathbb{E}[z(t+\tau)z(t)]$ do not depend on $t$~\cite[Ch.~9]{leongarcia08}.}
\begin{align*}
\mathbb{E}[y_i(t+\tau)y_i(t)]&= h_i(\tau)*h_i(-\tau)*\mathbb{E}[x_i(t+\tau)x_i(t)]\\
&= h_i(\tau)*h_i(-\tau)\\
&= \int_{-\infty}^{+\infty}h_i(-\nu) h_i(\tau-\nu) \d \nu.
\end{align*}
By setting $\tau=0$ and upon skipping some mundane integral evaluations, the variance of $y_i(t)$ for $i\geq 1$ can be computed as
\begin{equation}\label{eq:var1}
\mathbb{E}[y_i^2(t)]=\int_{0}^{\infty} h_i^2(\nu)\d \nu = \frac{1}{2\lambda_i \gamma}.
\end{equation}

We propose capturing the energy of inter-area oscillations through the metric
\begin{equation}\label{eq:metric1}
f_y:=\mathbb{E}[\|\by_\mcE\|_2^2]=\sum_{i\in\mcE}\mathbb{E}[y_i^2(t)]=\frac{1}{2\gamma}\trace(\bLambda_\mcE^{-1}).
\end{equation}
Ideally, one could consider the metric
\begin{equation*}
f_\delta:=\mathbb{E}[\|\bdelta_\mcE\|_2^2]=\frac{1}{2\gamma}\trace(\bM^{-1/2}\bU_\mcE\bLambda_\mcE^{-1} \bU_\mcE^\top\bM^{-1/2})
\end{equation*}
to account for the effect of transforming $\by_{\mcE}$ by $\bM^{-1/2}\bU_\mcE$ to obtain $\bdelta_\mcE$. However, incorporating $f_\delta$ in an optimization seems to be challenging. By definition $\bdelta_\mcE=\bM^{-1/2}\bU_\mcE\by_{\mcE}$ and the sub-multiplicative property of norms, metric $f_y$ can be shown to bound $f_\delta$ as
\[\frac{1}{\sqrt{M_\text{max}}} f_y\leq f_\delta\leq \frac{1}{\sqrt{M_\text{min}}} f_y\]
where $M_\text{max}$ and $M_\text{min}$ are the maximum and minimum inertia coefficients accordingly. 

We henceforth focus on dealing with $f_y$. The trace operator in \eqref{eq:metric1} returns the sum of the $K$ largest eigenvalues of $\bL_M^\dagger$. Matrix $\bL_M^\dagger$ in turn depends on the steady-state point of the power system. Hence, the system operator may consider selecting the steady-state point of the power system to minimize the stability metric $f_y$. To include $f_y$ into an OPF, the next result provides an alternative parameterization of $f_y(\bL_M)$, which is in fact a \emph{convex} function of $\bL_M$; see the Appendix for proof. 
\begin{lemma}\label{le:fy}
Given matrix $\bL_M$, the stability metric $f_y$ can be expressed as the optimal value of the ensuing SDP
\begin{subequations}\label{eq:sumEmax}
	\begin{align}
	f_y(\bL_M)=\min_{s,\bZ\succeq 0}\ &~\frac{\trace(\bZ)+Ks}{2\gamma}\label{eq:sumEmax:1}\\
	\mathrm{s.to}\ &~\begin{bmatrix}
	\bZ+s\bI_N & \bW\\
	\bW &\bL_M
	\end{bmatrix}\succeq 0\label{eq:sumEmax:2}
	\end{align} 
\end{subequations}
where $\bW:=\bI_N-\bu_1\bu_1^\top$ with $\bu_1$ defined in \eqref{eq:u1}.
\end{lemma}

We next show how this reformulation for $f_y(\bL_M)$ can be integrated into an SDP relaxation of the OPF.

\section{Inter-area Oscillations-Aware OPF}\label{sec:opf}
Given demands, this section devises an OPF to find an operating point balancing between the minimum generation cost and the stability-aware metric $f_y$. The complex power $p_n+jq_n$ at bus $n\in\mcN$ relates to external voltages $\bv$ as
	\begin{subequations}\label{eq:pfc}
		\begin{align}
		p_n&=\bv^{H}\bM_{p_n}\bv\label{eq:pfc:p}\\
		q_n&=\bv^{H}\bM_{q_n}\bv\label{eq:pfc:q}
		\end{align}
	\end{subequations}
where $(\bM_{p_n},\bM_{q_n})$ are $N\times N$ Hermitian matrices\footnote{Note the notational distinction from the matrix $\bM$ of inertia constants.}; see e.g.,~\cite{redux}. Squared voltage magnitudes and apparent line currents can be expressed as quadratic functions of $\bv$ too
\begin{align}\label{eq:volts}
V_n^2&=\bv^{H}\bM_{v_n}\bv\\
|\tilde{i}_{mn}|^2&=\bv^{H}\bM_{i_{mn}}\bv
\end{align}
where matrices $\bM_{v_n}$ and $\bM_{i_{mn}}$ are defined in~\cite{L2O}. 

The active power injected into bus $n$ can be decomposed into a dispatchable generation $p_n^g$ and fixed demand $p_n^d$ as $p_n=p_n^g-p_n^d$. Flexible demands at bus $n$, if any, can also be incorporated in $p_n^g$. The reactive power injected into bus $n$ is decomposed similarly as $q_n=q_n^g-q_n^d$. Recall that $p_n^g=p_n^d=q_n^g=q_n^d=0$ for all buses $n\in\bar{\mcS}$. Given the demands at all buses $\{p_n^d,q_n^d\}_{n\in\mcS}$, the OPF problem aims at optimally dispatching generators $\{p_n^g,q_n^g\}_{n\in\mcS}$ while meeting resource and network limits. 

The stability-aware OPF can be formulated as 
\begin{subequations}\label{eq:OPF1}
	\begin{align}
	\min\ &~(1-\mu)\left(\sum_{n\in\mcS} c_n^p p_n^g+c_n^q
	q_n^g\right)+\mu f_y(\bL_M)\nonumber\\
	\mathrm{over}\ &~\bv, \be, \bL, \{p_n^g,q_n^g\}_{n\in\mcS}\notag\\
	\mathrm{s.to}\ &~\eqref{eq:ve},~\eqref{eq:Laplacian},~\eqref{eq:L-LM}~\text{and}\nonumber\\
	&~\underline{p}_n^g\leq p_n^g\leq \bar{p}_n^g,\forall~n\in\mcS\label{eq:P1:pglim}\\
	&~\underline{q}_n^g\leq q_n^g\leq \bar{q}_n^g,\forall~n\in\mcS\label{eq:P1:qglim}\\
	&~\bv^{H}\bM_{p_n}\bv=p_n^g-p_n^d,\forall~n\in\mcS\label{eq:P1:pg}\\
	&~\bv^{H}\bM_{q_n}\bv=q_n^g-q_n^d,\forall~n\in\mcS\label{eq:P1:qg}\\
	&~\underline{v}_n\leq \bv^H\bM_{v_n}\bv\leq \bar{v}_n,\forall~n\in\mcN\label{eq:P1:vlim}\\	
	&~\bv^{H}\bM_{i_{mn}}\bv\leq \bar{i}_{mn},\forall~(m,n)\in\mcL\label{eq:P1:flim}
	\end{align} 
\end{subequations}
where $\mu\in[0,1]$ is a trade-off parameter, and $(c_n^p,c_n^q)$ are the coefficients for generation cost at bus $n$. Constraints \eqref{eq:P1:pglim}--\eqref{eq:P1:qglim} impose generation limits. Constraints \eqref{eq:P1:pg}--\eqref{eq:P1:qg} enforce the power flow equations at buses with non-zero injections. Constraints \eqref{eq:P1:vlim}-\eqref{eq:P1:flim} confine voltage magnitudes and apparent line currents within line ratings. 

The standard OPF formulation corresponds to $\mu=0$ and is known to be non-convex. It is also known however that this non-convex problem can be handled via its semidefinite program relaxation~\cite{Bai08}, which turns out to be exact under a broad range of operating conditions~\cite{Lavaei12}. The novel points in \eqref{eq:OPF1} are: \emph{i)} The stability-aware metric $f_y$, which by Lemma~\ref{le:fy} can be formulated as an SDP over the matrix variable $\bL_M$ or equivalently $\bL$; \emph{ii)} Matrix $\bL$ depends on internal voltages $\be$ via \eqref{eq:Laplacian}; and \emph{iii)} Internal voltages $\be$ have been included as variables and are linearly related to external voltages $\bv$ through \eqref{eq:ve}.

The crux of \eqref{eq:OPF1} seems to be the complicated nonlinear dependence of $\bL$ on the magnitudes and angles of $\be$. Surprisingly, this dependence can be dealt with an SDP relaxation: Introduce the \emph{lifted} matrix variable $\bE:=\be\be^H$. Thanks to $\bE$, constraint~\eqref{eq:Laplacian} is equivalent to
\begin{equation}\label{eq:Laplacian2}
[\bL]_{n,m}=\left\{\begin{array}{ll}
\sum_{k\neq n} \frac{\real\{E_{nk}\}}{\gamma_{nk}}&,~m=n\\
\\
-\frac{\real\{E_{nm}\}}{\gamma_{nm}}&,~ m\neq n.
\end{array}\right.
\end{equation}
Constraint~\eqref{eq:Laplacian2} relates $\bL$ to the lifted variable $\bE$ via linear equality constraints.

Using Lemma~\ref{le:fy} and the reformulation in~\eqref{eq:Laplacian2}, the SDP-relaxation of~\eqref{eq:OPF1} can be posed as
\begin{subequations}\label{eq:OPF2}
	\begin{align}
	\min\ &(1-\mu)\left(\sum_{n\in\mcS} c_n^p p_n^g+c_n^q
	q_n^g\right)+\mu\left(\frac{\trace(\bZ)+Es}{2\gamma}\right)\nonumber\\
	\mathrm{s.to}\ &~\eqref{eq:P1:pglim}-\eqref{eq:P1:qglim},~\eqref{eq:Laplacian2},~\text{and}\nonumber\\
	&~\trace(\bV\bM_{p_n})=p_n^g-p_n^d,\forall~n\in\mcS\label{eq:SR:pg}\\
	&~\trace(\bV\bM_{q_n})=q_n^g-q_n^d,\forall~n\in\mcS\label{eq:SR:qg}\\
	&~\underline{v}_n\leq \trace(\bV\bM_{v_n})\leq \bar{v}_n,\forall~n\in\mcN\label{eq:SR:vlim}\\	
	&~\trace(\bV\bM_{i_{mn}})\leq \bar{i}_{mn},\forall~(m,n)\in\mcL\label{eq:SR:flim}\\
	&~\bV_{\mcS,\mcS}=\bGamma\bY_{\mcS}\bE\bY_{\mcS}^H\bGamma^H,\label{eq:SR:VE}\\
	&~\begin{bmatrix}
	\bZ+s\bI_N & \bW\\
	\bW &\bM^{-1/2}\bL\bM^{-1/2}
	\end{bmatrix}\succeq 0,\label{eq:SR:epigraph}\\
	&~\bV\succeq 0\label{eq:SR:Vpsd}
	\end{align} 
\end{subequations}
where $\bV_{\mcS,\mcS}$ is the submatrix of $\bV$ obtained from its first $S$ rows and columns. Constraint \eqref{eq:SR:Vpsd} implies $\bV_{\mcS,\mcS}\succeq 0$, which combined with \eqref{eq:SR:VE} implies also that $\bE\succeq 0$. As in the standard SDP of the OPF, problem \eqref{eq:OPF2} is termed \emph{exact} if the optimal $\bV$ and consequently $\bE$ are rank-1. Our numerical tests demonstrate that \eqref{eq:OPF2} is exact for all instances tested.

\section{Numerical Tests}\label{sec:tests}
The proposed approach for generation dispatch while minimizing low-frequency rotor angle oscillations was tested on the IEEE 39-bus system. All OPF instances were solved using the MATLAB-based optimization toolbox
YALMIP alongside the SDP solver Mosek~\cite{YALMIP}, \cite{mosek}. Network parameters and nominal loads were taken from the MATPOWER casefile~\cite{MATPOWER}. The IEEE 39-bus system features 10 generators and 10 zero-injection buses. After eliminating the zero-injection buses, a Kron-reduced system with $|\mcS|=29$ buses was obtained. Inertia coefficients for generators were taken from~\cite{canizares2015benchmark} and for the remaining nonzero-injection buses in $\mcS$ were set to $10\%$ of the average generator inertia. The proportionality coefficient $\gamma$ was set to $0.1467$ based on the ratio of average damping to average inertia in~\cite{BDNK19}. Since the benchmark system has identical generation cost curves across generators, a uniform (re)active power cost $c_n^p=1$ and $c_n^q=0.1$  was used for all generators. Reactive power costs were included to ensure the exactness of the SDP relaxation~\cite{MSL15}.

For $\mu=0$, problem~\eqref{eq:OPF2} solves the conventional OPF to minimize generation cost. To verify if the relaxation remains exact for $\mu>0$, we evaluated the ratio of the second largest to the largest eigenvalue of $\bV$ for all tests. We found these ratios to be less than $10^{-3}$ for all test instances. This indicates that~\eqref{eq:OPF2} features nearly rank-1 $\bV$, thus rendering the relaxation exact.

\begin{figure}[thb]
    \centering
	\includegraphics[width=1.\linewidth]{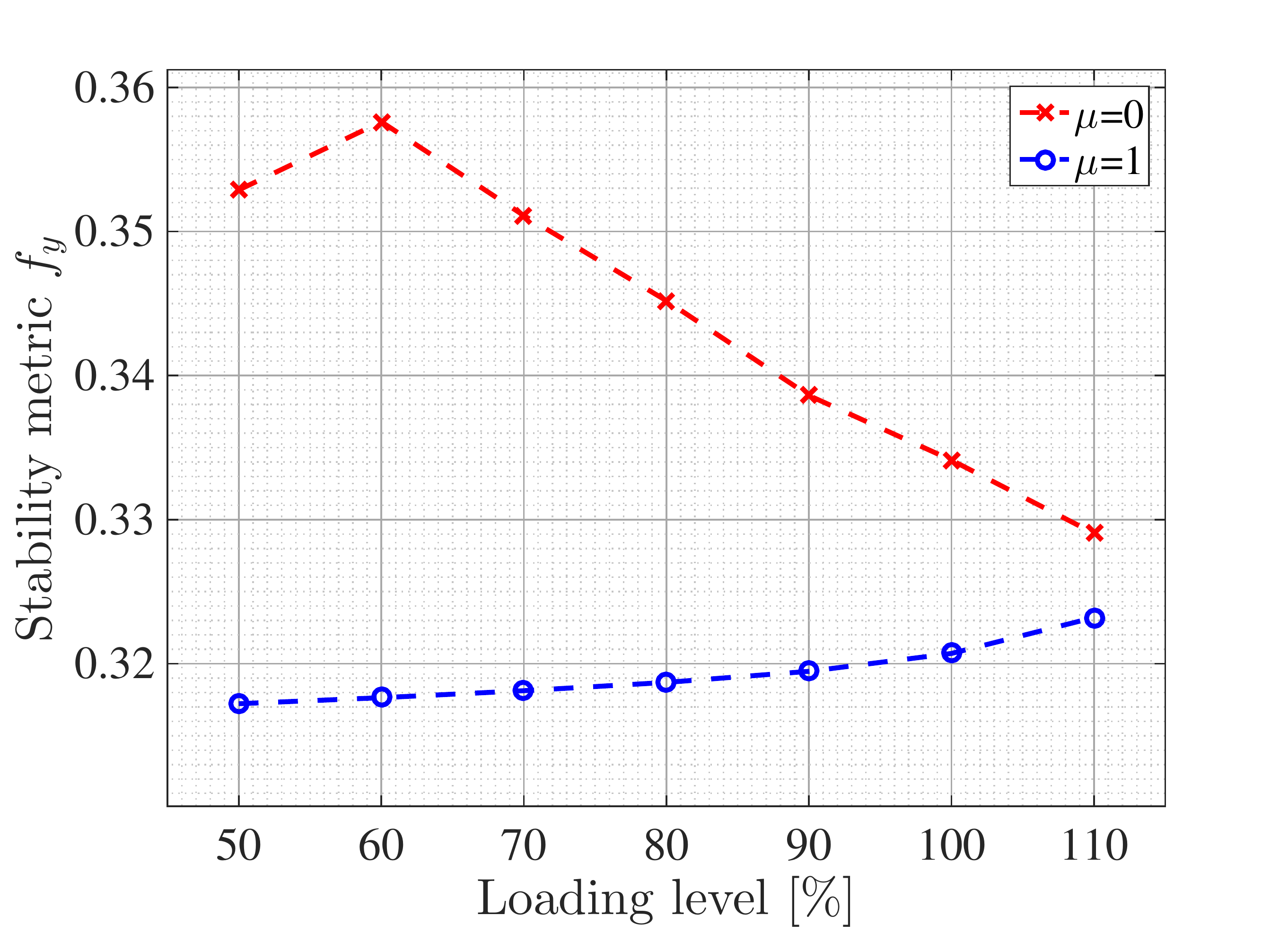}
	\caption{Stability metric $f_y$ for $K=3$ eigensystems for increasing load levels attained by a cost-focused ($\mu=0$) and a stability-focused ($\mu=1$) OPF.}
	\label{fig:loadlevel}
\end{figure}

\begin{table}[thb]
\centering
\caption{\label{font-table} Stability metric $f_y$ for different number $K$ of eigensystems at $50\%$ loading. \vskip 3pt }
\label{tab:table1}
\begin{tabular}{l|l|P{2.5cm}}
\hline
\multicolumn{2}{l|}{} & $f_y~[\times 10^{-2}]$  \\ \hline
\multirow{2}{*}{$K=1$} & $\mu=0$  & 19.19  \\ \cline{2-3} 
                   & $\mu=1$  & 18.03   \\ \hline
\multirow{2}{*}{$K=2$}  & $\mu=0$  & 27.81  \\ \cline{2-3} 
                   & $\mu=1$  & 25.42  \\ \hline
\multirow{2}{*}{$K=3$}  & $\mu=0$ & 35.30  \\ \cline{2-3} 
                   & $\mu=1$  & 31.73  \\ \hline
\multirow{2}{*}{$K=4$}  &  $\mu=0$ & 41.88  \\ \cline{2-3} 
                   & $\mu=1$  & 37.28  \\ \hline
\multirow{2}{*}{$K=5$}  &  $\mu=0$ & 46.96  \\ \cline{2-3} 
                   & $\mu=1$  & 41.85   \\ \hline
\end{tabular}
\end{table}

We examined the improvement in terms of the stability metric $f_y$ in~\eqref{eq:metric1} obtained while transitioning from the traditional OPF that minimizes generation cost to the proposed stability-improving OPF. We first considered the first $K=3$ eigensystems with the lowest frequencies. Figure~\ref{fig:loadlevel} plots $f_y$ for $\mu=0$ and $\mu=1$, and by scaling (re)active loads from $50\%$ to $110\%$ of nominal values. From these test instances, it was found that the proposed approach improves $f_y$ by up to $11.18\%$ over the traditional approach. Another interesting observation from Figure~\ref{fig:loadlevel} is that improvements in stability are more prominent at lower loading conditions. The previous analysis can be extended to any number $K$ of low-frequency eigensystems one seeks to consider. To illustrate that, Table~\ref{tab:table1} reports the values of $f_y$ attained for different values of $K$, all at $50\%$ load level.

\begin{figure}[thb]
    \centering
	\includegraphics[width=1.\linewidth]{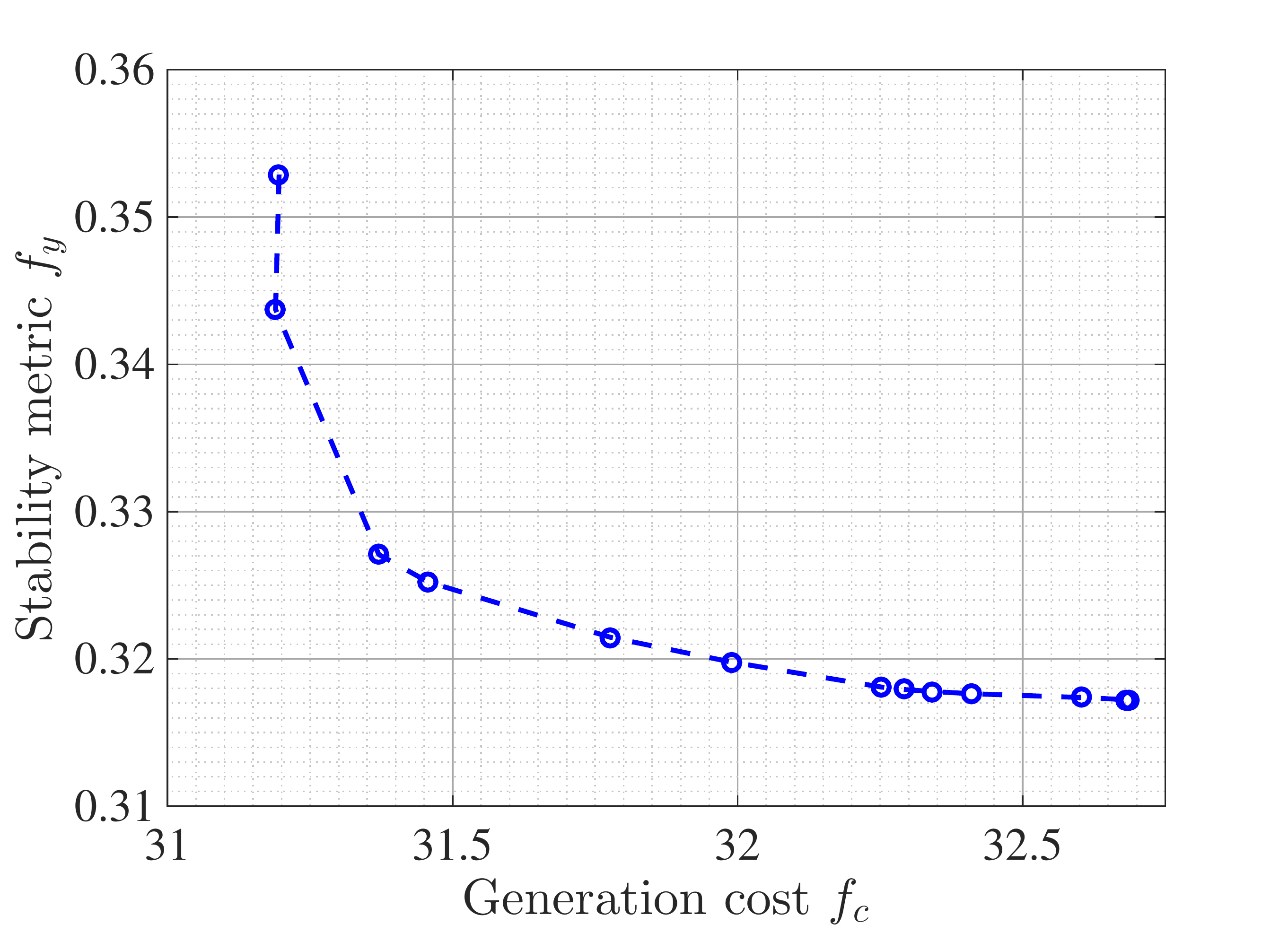}
	\caption{Generation cost and stability trade-off at $50\%$ loading and $K=3$.}
	\label{fig:pareto}
\end{figure}

Having gained an insight in the overall room for improvement in stability for different load levels and $K$, we then focused on the trade-off between stability and generation cost for a fixed loading setup. We conducted a Pareto optimality analysis by varying $\mu\in[0,1]$ with (re)active power loads fixed at $50\%$ of their nominal value. Figure~\ref{fig:pareto} shows the stability metric $f_y$ for $K=3$ versus the total (re)active power generation cost $f_c:=\sum_{n\in\mcS} c_n^p p_n^g+c_n^q q_n^g$. The Pareto analysis shows that an improvement of $10.14\%$ in $f_y$ can be achieved by increasing the generation cost by $4.78\%$.
However, instead of seeking the most stable operating point associated with $\mu=1$, one can choose an intermediate operating point on the Pareto front to significantly improve stability with marginal increase in generation cost. For instance, the stability metric $f_y$ can be reduced by $7.8\%$ by increasing the generation cost only by $0.84\%$. Therefore, the proposed approach provides a flexible tool for power grid dispatchers to include stability considerations based on willingness to spend on additional generation cost.

\section{Conclusions}\label{sec:conclusions}
A detailed model capturing the dependence of power system small-signal dynamics on the operating point has been put forth. Thereafter, a novel OPF formulation for minimizing the energy of low-frequency inter-area oscillations has been proposed. To garner global optimality and scalability, the proposed OPF relies on SDP relaxation, which can be conveniently integrated with classical OPF formulations that minimize generation costs. Thanks to such integration, an insightful investigation of the trade-off between the chosen stability metric and generation cost was carried out for the IEEE 39-bus benchmark system. The numerical tests interestingly reveal that despite relying on SDP relaxation, the proposed stability-improving OPF yields nearly rank-1 minimizers, rendering the relaxation exact. The numerical tests further exemplify that the developed OPF can be used to obtain an optimally stable operating point based on the number of low-frequency modes one intends to suppress and the acceptable deviation from least generation cost dispatch. This work sets the foundations for interesting research directions. Improving small-signal stability via reactive power dispatch is a lucrative next step, mainly due to its limited monetary requirements. Expanding the proposed framework to include diverse stability metrics and more detailed generator models would further help to expand its applicability.

\appendix
\begin{IEEEproof}[Proof of Lemma~\ref{le:fy}]
It is known that the sum of the $K$ largest eigenvalues of matrix $\bX$ can be expressed as the SDP~\cite[pp.~54]{VaBo96}
\begin{align*}
\min_{s,\bZ\succeq 0}\ &~\trace(\bZ)+Ks\\
\mathrm{s.to}\ &~\bZ+s\bI_N\succeq \bX.
\end{align*} 
To show Lemma~\ref{le:fy}, it suffices to show that constraint \eqref{eq:sumEmax:2} is equivalent to 
\begin{equation}\label{eq:LMI1}
\bZ+s\bI_N\succeq \bL_M^\dagger.
\end{equation}
This follows from the Schur complement of the left-hand side of \eqref{eq:sumEmax:2} with respect to the singular block $\bL_M$~\cite[Sec.~A.5.5]{BoVa04}. More specifically, constraint \eqref{eq:sumEmax:2} is equivalent to the three conditions
\begin{align*}
&\bL_M\succeq 0\\
&(\bI_N-\bL_M\bL_M^\dagger)\bW=\bzero\\
&\bZ+s\bI_N\succeq \bW\bL_M^\dagger\bW.
\end{align*}
The first condition holds by definition of $\bL_M$ and for practical voltage angle differences across neighboring buses. 

For the second condition, let $\bU_2$ be the submatrix of $\bU$ collecting the eigenvectors associated with the non-zero eigenvalues of $\bL_M$. If the non-zero eigenvalues are placed on the main diagonal of matrix $\bLambda_2$, the eigenvalue decomposition yields $\bL_M=\bU_2\bLambda_2\bU_2^\top$, $\bL_M^\dagger=\bU_2\bLambda_2^{-1}\bU_2^\top$, and $\bL_M\bL_M^\dagger=\bU_2\bU_2^\top$.

The orthonormality of eigenvectors provides that $\bu_1\bu_1^\top+\bU_2\bU_2^\top=\bI_N$, so $\bW$ is expressed as 
\begin{equation}\label{eq:W2}
\bW=\bI_N-\bu_1\bu_1^\top=\bU_2\bU_2^\top.
\end{equation}
Then, the second condition is satisfied as $(\bI_N-\bU_2\bU_2^\top)\bW=\bzero$. Because of \eqref{eq:W2}, we obtain $\bW\bL_M^\dagger\bW=\bL_M^\dagger$, so the third condition yields \eqref{eq:LMI1}.
\end{IEEEproof}

\balance

\bibliographystyle{ieeetr}
\bibliography{sample,myabrv,power}
\end{document}